\documentclass[12pt]{amsart}
\usepackage{amsfonts, amssymb, amsmath, amsthm, euscript}
\swapnumbers
\theoremstyle{plain}
\newtheorem{thm}{Theorem}[section]
\newtheorem{lem}[thm]{Lemma}
\newtheorem{prop}[thm]{Proposition}
\newtheorem{cor}[thm]{Corollary}

\theoremstyle{definition}
\newtheorem{setup}[thm]{Setup}
\newtheorem{exam}[thm]{Example}
\newtheorem{specialcase}[thm]{A special case}
\newtheorem{catenarity}[thm]{Catenarity}
\newtheorem*{Ack}{Acknowledgements}

\theoremstyle{remark}

\newtheorem{note}[thm]{}

\def\G{\Gamma}

\def\Cdim{\operatorname{Cdim}}
\def\rank{\operatorname{rank}}
\def\spec{\operatorname{Spec}}
\def\gr{\operatorname{gr}}

\def\qij{q_{ij}}
\def\nxn{n{\times}n}

\def\N{\mathbb{N}}
\def\Nn{\mathbb{N}^n}
\def\Z{\mathbb{Z}}

\def\Zn{\mathbb{Z}^n}
\def\F{\mathbb{F}}
\def\C{\mathcal{C}}
\def\kq{k_q}
\def\qpsx{k_q[[x]]}
\def\qlsx{k_q[[x^{\pm 1}]]}
\def\qpx{k_q[x]}
\def\qlx{k_q[x^{\pm 1}]}

\def\ktimes{k^{\times}}
\def\ltgrlex{\; {<}_\text{grlex} \; }

\def\tildeA{\widetilde{A}}

\begin{document}

\title{Prime Ideals of $q$-Commutative Power Series Rings}

\author{Edward S. Letzter}

\author{Linhong Wang}

\address{Department of Mathematics\\
  Temple University\\
  Philadelphia, PA 19122}

\email{letzter@temple.edu}

\address{Department of Mathematics\\
  Southeastern Louisiana University\\
  Hammond, LA 70402}

\email{lwang@selu.edu}

      \thanks{Research of the first author supported in part by grants
        from the National Security Agency. Results included in this
        paper appear in the second author's
        Ph.D. thesis at Temple University.}



\begin{abstract} We study the ``$q$-commutative'' power series ring
  $R:=\kq[[x_1,\ldots,x_n]]$, defined by the relations $x_ix_j = \qij
  x_j x_i$, for mulitiplicatively antisymmetric scalars $\qij$ in a
  field $k$. Our results provide a detailed account of prime
  ideal structure for a class of noncommutative, complete, local,
  noetherian domains having arbitrarily high (but finite) Krull,
  global, and classical Krull dimension. In particular, we prove that
  the prime spectrum of $R$ is normally separated and is finitely
  stratified by commutative noetherian spectra. Combining this normal
  separation with results of Chan, Wu, Yekutieli, and Zhang, we are
  able to conclude that $R$ is catenary. Following the approach of
  Brown and Goodearl, we also show that links between prime ideals are
  provided by canonical automorphisms. Moreover, for sufficiently
  generic $\qij$, we find that $R$ has only finitely many prime ideals
  and is a UFD (in the sense of Chatters).
\end{abstract}

\maketitle

%%%%%%%%%%%%%%%%%%%%%%%%%%%%%%%%%%-Begin-Body-%%%%%%%%%%%%%%%%%%%%%%%%%%%%%%%%%

\section{Introduction}

Given a field $k$ and an $\nxn$ matrix $q = (\qij)$, with $q_{ii} = 1$
and $\qij = q^{-1}_{ji} \in \ktimes$, we can construct the
``$q$-commutative'' power series ring $R =
\kq[[x_1,\ldots,x_n]]$. Multiplication is determined by the
commutation relations $x_ix_j = \qij x_j x_i$, leading to a rich
combinatorial structure (see, e.g., \cite{Koo}). Moreover, it follows
from longstanding ring-theoretic results and techniques
(cf.~\cite{LiVOy},\cite{Wal}) that $R$ is a complete, local, Auslander
regular, noetherian, zariskian (in the sense of Li and Van Oystaeyen
\cite{LiVOy}) domain with Krull dimension, classical Krull dimension,
and global dimension all equal to $n$. In this paper, the goal is to
provide a detailed account of the prime ideal theory of $R$. Our
approach builds on earlier studies (particularly \cite[\S 4]{BroGoo2},
\cite[\S 2]{DeCKacPro}, and \cite{GooLet}) on $q$-commutative
polynomial rings (also known as ``quantum affine spaces'' and
``twisted polynomial algebras''). However, it should be no surprise
that our work in this paper frequently involves topological
considerations not needed in the earlier studies.

\begin{note}
  Our results depend on a careful examination of the two-sided ideal
  structure of the $q$-commutative Laurent series ring $L=\kq[[x^{\pm
    1}_1,\ldots,x^{\pm 1}_n]]$. This analysis can be
  briefly described as follows: \smallskip

  To start, there is an obvious action of the $n$-torus $H=
  (\ktimes)^n$ on $L$ (and $R$) by automorphisms. In (\ref{H-simple})
  we prove, when $k$ is infinite, that $L$ is $H$-simple (i.e., that
  the only $H$-stable ideals of $L$ are the zero ideal and $L$
  itself). Consequently, every $H$-orbit of prime ideals of $L$ is
  Zariski dense in $\spec L$ (when $k$ is infinite).  \smallskip

  Next, in (\ref{centrally_generated}) we prove that extension and
  contraction of ideals produces a bijection
%%%%%%%%%%%%%%%%%%%%%%%%%%%%%%%%%%%%%%%%%%%%%%%%%%%%%%%%%%%%%%%%%%%%%%%%%%%%%%%
  \[\left\{ \text{ideals of $L$} \right\} \longleftrightarrow
  \left\{\text{ideals of the center $Z = Z(L)$} \right\} \]
%%%%%%%%%%%%%%%%%%%%%%%%%%%%%%%%%%%%%%%%%%%%%%%%%%%%%%%%%%%%%%%%%%%%%%%%%%%%%%%
  for any choice of the field $k$. This bijection produces a
  homeomorphism between the prime spectrum $\spec L$, equipped with
  the Zariski topology, and $\spec Z$; see (\ref{homeomorphism}). It
  follows that $Z$ is a noetherian domain.  \smallskip

  Our results for $L$ parallel, in part, those found in \cite[\S
  4]{BroGoo2}, \cite[\S 2]{DeCKacPro}, and \cite[\S 1]{GooLet} for
  $q$-commutative Laurent polynomial rings (also known as ``quantum
  tori'' and ``McConnell-Pettit algebras''). Additional references are
  given in (\ref{lemma_history}). However, in contrast to our own work
  below, the $q$-commutative Laurent polynomial case permits changes
  of variables not available for $q$-commutative Laurent series; see
  (\ref{change_of_variables}). In particular, while the center of a
  $q$-commutative Laurent polynomial ring is isomorphic to a
  commutative Laurent polynomial ring, it is possible (following an
  observation of K. R. Goodearl) that $Z$ as above is not a
  commutative Laurent series ring; see (\ref{not_Laurent}).
\end{note}

\begin{note} The analysis of $L$ can be applied to $R$ as
  follows. First of all, the $x_1,\ldots,x_n$ provide a stratification
%%%%%%%%%%%%%%%%%%%%%%%%%%%%%%%%%%%%%%%%%%%%%%%%%%%%%%%%%%%%%%%%%%%%%%%%%%%%%%%
\[ \spec R  =  \bigsqcup _{w \in W} \spec_wR ,\]
%%%%%%%%%%%%%%%%%%%%%%%%%%%%%%%%%%%%%%%%%%%%%%%%%%%%%%%%%%%%%%%%%%%%%%%%%%%%%%%
where each $w$ is a subset of $\{1,\ldots,n\}$, and where
%%%%%%%%%%%%%%%%%%%%%%%%%%%%%%%%%%%%%%%%%%%%%%%%%%%%%%%%%%%%%%%%%%%%%%%%%%%%%%
\[ \spec_w R  =  \left\{ P \in \spec R \mid x_i \in P
  \Leftrightarrow i \in w \right\}.\]
%%%%%%%%%%%%%%%%%%%%%%%%%%%%%%%%%%%%%%%%%%%%%%%%%%%%%%%%%%%%%%%%%%%%%%%%%%%%%%%
Each $\spec_w R$ is naturally homeomorphic to $\spec L_w$, where $L_w$
is a $q$-commutative Laurent series ring in $n-\vert w \vert$
variables, for a suitable replacement of the original matrix $q$; see
(\ref{stratification}).
\smallskip

Furthermore, each $\spec _wR$ is a union of $H$-orbits in $\spec R$,
and every $H$-orbit in $\spec_wR$ is dense in $\spec_wR$ (when $k$ is
infinite); see (\ref{H-equivariance}).
We also see, in (\ref{Jw}), that the $H$-prime ideals of $R$ are
exactly those ideals of the form $\langle x_i \mid i \in w \rangle$,
for subsets $w$ of $\{1,\ldots,n\}$.
\smallskip

This stratification of $\spec R$ parallels that found in \cite[\S
4]{BroGoo2}, \cite[\S 2]{DeCKacPro}, and \cite[\S 2]{GooLet} for
$q$-commutative polynomial rings.
\end{note}

\begin{note} Noetherian unique factorization domains were defined by
  Chatters \cite{Cha} (cf.~\cite{ChaJor}) to be noetherian domains
  such that each height-one prime ideal is completely prime and
  generated by a normal element. In \cite{Ven}, Venjakob presented the
  first known examples of noetherian, noncommutative, complete, local,
  unique factorization domains of global, Krull, and classical Krull
  dimension $d > 1$; these examples are completed group algebras of
  certain uniform pro-$p$ groups of rank $2$, and for these examples
  $d = 2$. 

  In (\ref{special_case}) we note, for sufficiently generic choices of
  the $\qij$, that the prime ideals of $R$ are precisely the $H$-prime
  ideals mentioned above: $\langle x_i \mid i \in w \rangle$,
  for subsets $w$ of $\{1,\ldots,n\}$. It then follows, in this
  special case, that $R$ is a noetherian unique factorization domain,
  in Chatters' sense, having global, Krull, and classical Krull
  dimension $n$.
\end{note} 

\begin{note} Following the approach of Brown and Goodearl in
  \cite{BroGoo2}, we also consider issues related to the localization
  and representation theory of $R$ (see \cite{BroWar} and
  \cite{Jat}). First, in (\ref{normal_and_link}i) we show that $\spec
  R$ is normally separated. Consequently, $R$ satisfies the strong
  second layer condition. Second, letting $G$ denote the group of
  automorphisms of $R$ generated by $r \mapsto x_i r x^{-1}_i$, for $r
  \in R$ and $1 \leq i \leq n$, we show in (\ref{normal_and_link}ii)
  that if $P \rightsquigarrow Q$ in $\spec R$ then there exists an
  automorphism $\tau \in G$ such that $\tau(P) = Q$. It then follows
  from well-known results of Stafford \cite{Sta} and Warfield
  \cite{War} that $R$ has a classical localization theory, in the
  sense of Jategaonkar \cite{Jat}, when $k$ is uncountable.
\end{note}

\begin{note} Combining the normal separation for $R$ established in (\ref{normal_and_link}i) with results of Chan, Wu, Yekutieli, and Zhang \cite{ChaWuZha}, \cite{WuZha}, \cite{YekZha}, we are able to conclude that $R$ is catenary.
\end{note}

\begin{note}
  In addition to similarities with $q$-commutative polynomial rings,
  the behavior of $q$-commutative power series rings also has deep
  analogues in the theory of quantum semisimple groups and other
  quantum function algebras; see (e.g.)
  \cite{BroGoo1},\cite{BroGoo2},\cite{HodLevTor},\cite{Jos},
  \cite{LauLen}.
\end{note}

\begin{note} As noted above, $R$ is zariskian (with respect to the
  $\langle x_1,\ldots,x_n\rangle$-adic filtration), and so it follows
  from \cite{LiVOy} that the ideals of $R$ are closed in the $\langle
  x_1,\ldots,x_n\rangle$-adic topology; see (\ref{I-adic}ii) and
  (\ref{closed}). This fact plays a key role throughout this paper.
\end{note}

\begin{note} The Goldie ranks of the prime factors of a given
  $q$-commutative polynomial ring have a single finite upper bound
  (see \cite{GooLet}). At this time we can only ask whether a similar
  upper bound exists, in general, for the Goldie ranks of the prime
  factors of $R$. 
\end{note}

\begin{note} The paper is organized as follows: Section 2 reviews some
  well known results and develops some of the preliminary
  theory. Section 3 is a thorough analysis of the ideal structure of
  $L$. Section 4 gives the main results for $R$. 
\end{note}

\begin{Ack} We are grateful to James Zhang for pointing out to us
  Corollary \ref{cat} on the catenarity of $q$-commutive power series
  rings, and we are also grateful to him for suggesting to us the
  proof provided for this corollary. We are grateful to Ken
  Goodearl for his observation noted in (\ref{not_Laurent}).
\end{Ack}

\section{$q$-Commutative Power and Laurent Series Rings:
Preliminaries}

We assume basic familiarity with noetherian rings (see, e.g.,
\cite{GooWar},\cite{McCRob}) and filtered rings (see, e.g.,
\cite{LiVOy}, \cite{NasVOy}).

\begin{setup} \label{setup} The following notation will remain in effect
  throughout this paper.

  (i) To start, $k$ will denote a field, $\ktimes$ will denote the
  multiplicative group of units in $k$, $n$ will denote a positive
  integer, and $q = (\qij)$ will denote a multiplicatively
  antisymmetric $\nxn$ matrix (i.e., $q_{ii} = 1$ and $\qij =
  q^{-1}_{ji}$). We will use $\N$ to denote the set of non-negative
  integers. We will use $\spec$ to denote the set of prime ideals of a
  designated ring, and we will always equip such sets with the Zariski
  topology.

  (ii) $R := \kq[[x_1,\ldots,x_n]]$ will denote the associative unital
  $k$-algebra of formal skew power series in the indeterminates
  $x_1,\ldots,x_n$, subject only to the commutation relations $x_i x_j
  = \qij x_j x_i$. We will refer to $R$, in general, as a
  \emph{$q$-commutative power series ring\/} (allowing $q$ and $n$ to
  vary in this usage).

  The elements of $R$ are power series
%%%%%%%%%%%%%%%%%%%%%%%%%%%%%%%%%%%%%%%%%%%%%%%%%%%%%%%%%%%%%%%%%%%%%%%%%%%%%%%
 \[ \sum _{s \in \Nn} c_s x^s ,\]
%%%%%%%%%%%%%%%%%%%%%%%%%%%%%%%%%%%%%%%%%%%%%%%%%%%%%%%%%%%%%%%%%%%%%%%%%%%%%%%
 for $c_s \in k$, for $s = (s_1,\ldots,s_n) \in \Nn$, and for $x^s =
 x_1^{s_1}\cdots x_n^{s_n}$. We will use $x^s$ to refer to a general
 monic monomial in $R$, and we can write $R = \qpsx$.

 (iii) Contained within $R$ is the $k$-subalgebra $\qpx =
 \kq[x_1,\ldots,x_n]$ of $q$-commuting polynomials in
 $x_1,\ldots,x_n$. This algebra has been extensively studied; see, for
 example, \cite[\S 4]{BroGoo2}, \cite[\S 2]{DeCKacPro},
 \cite{GooLet}.  As noted above, our analysis builds on these
 studies. Observe that $R$ is the completion of $\qpx$ at the ideal
 generated by $x_1,\ldots,x_n$.

 (iv) $J$ will denote the augmentation ideal $\langle x_1,\ldots , x_n
 \rangle$ of $R$.

 (v) $L$ will denote the $q$-commutative Laurent series ring $\qlsx =
 \kq[[x^{\pm 1}_1,\ldots,x^{\pm 1}_n]]$, defined later in
 (\ref{laurent}).
\end{setup}

Of course, when $k = \F_2$ it follows immediately that $R$ is the
commutative power series ring in $n$ variables. However, we will allow
this case unless indicated otherwise. We will also include the
possibility that $n=1$, in which case $R$ is the commutative power
series ring over $k$ in one variable.

\begin{note} \label{filtration} Much of our analysis depends on
  standard results concerning filtered rings, which we now
  briefly review. Again, the reader is referred to \cite{LiVOy} and
  \cite{NasVOy} for background.

  (i) To start, let $A$ be an (associative unital) ring, and 
  suppose further there exist additive
  subgroups
%%%%%%%%%%%%%%%%%%%%%%%%%%%%%%%%%%%%%%%%%%%%%%%%%%%%%%%%%%%%%%%%%%%%%%%%%%%%%%%
  \[ \cdots \supseteq A_{-2} \supseteq A_{-1} \supseteq A_0 \supseteq
  A_1 \supseteq A_2 \supseteq \cdots \]
%%%%%%%%%%%%%%%%%%%%%%%%%%%%%%%%%%%%%%%%%%%%%%%%%%%%%%%%%%%%%%%%%%%%%%%%%%%%%%%
  of $A$ with $A_iA_j \subseteq A_{i+j}$, for all integers $i$ and
  $j$. We refer to the preceding as a \emph{filtration\/} of $A$.
This filtration turns $A$ into a topological additive
  group by letting the cosets of the $A_i$, for all $i$, form a
  fundamental system of neighborhoods in $A$.

  (ii) The above filtration on $A$ is \emph{exhaustive\/} if
%%%%%%%%%%%%%%%%%%%%%%%%%%%%%%%%%%%%%%%%%%%%%%%%%%%%%%%%%%%%%%%%%%%%%%%%%%%%%%%
\[ A = \bigcup_{j \in \Z} A_j, \]
%%%%%%%%%%%%%%%%%%%%%%%%%%%%%%%%%%%%%%%%%%%%%%%%%%%%%%%%%%%%%%%%%%%%%%%%%%%%%%%%
is \emph{separated\/} if the intersection of the $A_j$ is
  equal to zero (or equivalently, if the corresponding topology is
  Hausdorff), and is \emph{complete\/} if Cauchy sequences converge in
  the corresponding topology.

  (iii) We have 
%%%%%%%%%%%%%%%%%%%%%%%%%%%%%%%%%%%%%%%%%%%%%%%%%%%%%%%%%%%%%%%%%%%%%%%%%%%%%%
  \[ \gr A = \cdots \oplus A_{-2}/A_{-1} \oplus A_{-1}/A_0 \oplus A_0/A_1 \oplus
  A_1/A_2 \oplus \cdots, \]
%%%%%%%%%%%%%%%%%%%%%%%%%%%%%%%%%%%%%%%%%%%%%%%%%%%%%%%%%%%%%%%%%%%%%%%%%%%%%%
the \emph{associated graded ring\/} corresponding to the given filtration.

(iv) Suppose that the filtration on $A$ is exhaustive, separated, and
complete. It follows from standard arguments that $A$ is right
noetherian if $\gr A$ is right noetherian and that $A$ is left noetherian
if $\gr A$ is left noetherian. 
\end{note}

\begin{note} \label{I-adic} Let $A$ be a ring, and let $I$ be an ideal of
  $A$.

  (i) Keeping (\ref{filtration}) in mind, we have the \emph{$I$-adic
    filtration\/} of $A$:
%%%%%%%%%%%%%%%%%%%%%%%%%%%%%%%%%%%%%%%%%%%%%%%%%%%%%%%%%%%%%%%%%%%%%%%%%%%%%%
\[ A = I^0 \supseteq I^1 \supseteq \cdots \]
%%%%%%%%%%%%%%%%%%%%%%%%%%%%%%%%%%%%%%%%%%%%%%%%%%%%%%%%%%%%%%%%%%%%%%%%%%%%%%
We also have the corresponding \emph{$I$-adic topology\/}, making $A$
a topological ring. Note that the $I$-adic filtration of $A$ is
exhaustive.

(ii) We also have the \emph{Rees ring\/} associated to the $I$-adic
filtration on $A$:
%%%%%%%%%%%%%%%%%%%%%%%%%%%%%%%%%%%%%%%%%%%%%%%%%%%%%%%%%%%%%%%%%%%%%%%%%%%%%%
\[ \tildeA := A \oplus I \oplus I^2 \oplus \cdots \]
%%%%%%%%%%%%%%%%%%%%%%%%%%%%%%%%%%%%%%%%%%%%%%%%%%%%%%%%%%%%%%%%%%%%%%%%%%%%%%
Following the definition in \cite[p.~83]{LiVOy}, when $I$ is contained
in the Jacobson radical of $A$, and when $\tildeA$ is right
noetherian, we say that $A$ is \emph{(right) zariskian\/} with respect
to the $I$-adic filtration.

In \cite[p.~87, Proposition]{LiVOy} it is proved that $A$ is (right)
zariskian with respect to the $I$-adic filtration if the $I$-adic
filtration is complete and if the associated graded ring corresponding
to the $I$-adic filtration is right noetherian.

The only consequence of the zariskian property necessary for our
analysis is the following: If $A$ is zariskian with respect to the
$I$-adic filtration, then every right ideal of $A$ is closed in the
$I$-adic topology \cite[p.~85, Corollary]{LiVOy}.
\end{note}

Before turning to $R$, we first consider skew power series, in a
single variable, over more general coefficient rings.

\begin{note} \label{one_variable} Let $A$ be a ring, let $\alpha$ be
  an automorphism of $A$, and let $A[[y;\alpha]]$ denote the ring of
  skew power series
%%%%%%%%%%%%%%%%%%%%%%%%%%%%%%%%%%%%%%%%%%%%%%%%%%%%%%%%%%%%%%%%%%%%%%%%%%%%%%%
  \[ \sum_{i = 0}^\infty a_i y^i  =  a_0 + a_1y + a_2y^2 + \cdots,\]
%%%%%%%%%%%%%%%%%%%%%%%%%%%%%%%%%%%%%%%%%%%%%%%%%%%%%%%%%%%%%%%%%%%%%%%%%%%%%%%
  for $a_0,a_1,\ldots \in A$, with multiplication determined by $ya =
  \alpha(a)y$, for $a \in A$. Of course, since $\alpha$ is an
  automorphism, we can just as well write the coefficients on the
  right. Set $B = A[[y;\alpha]]$.  As either a left or right
  $A$-module, we can view $B$ as a direct product of copies of $A$,
  indexed by $\N$.  Note that $y$ is normal in $B$ (i.e., $By = yB$),
  and $B/\langle y \rangle$ is naturally isomorphic to $A$.

  (i) Set $f = 1 + b_1y + b_2y^2 + \cdots \in B$, for $b_1,b_2,\ldots
  \in A$, and set $g = 1 + c_1y + c_2y^2 + \cdots \in B$, for
  $c_1,c_2,\ldots \in A$. Then
%%%%%%%%%%%%%%%%%%%%%%%%%%%%%%%%%%%%%%%%%%%%%%%%%%%%%%%%%%%%%%%%%%%%%%%%%%%%%%%
  \[ fg  =  1 + (b_1 + c_1)y + (b_2 + c_2 + p_2(b_1,c_1))y^2 +
  (b_3 + c_3 + p_3(b_1,b_2,c_1,c_2))y^3 + \cdots ,\]
%%%%%%%%%%%%%%%%%%%%%%%%%%%%%%%%%%%%%%%%%%%%%%%%%%%%%%%%%%%%%%%%%%%%%%%%%%%%%%%
  where $p_i(b_1,\ldots,b_{i-1},c_1,\ldots,c_{i-1}) \in A$ depends
  only on $b_1,\ldots,b_{i-1}$ and $c_1,\ldots,c_{i-1}$, for integers
  $i \geq 2$.  If $b_1,b_2,\ldots$ are arbitrary then we can choose
  $c_1,c_2,\ldots$ such that $fg = 1$, and if $c_1,c_2,\ldots$ are
  arbitrary then we can choose $b_1,b_2,\ldots$ such that $fg = 1$.

  We can deduce that a power series in $B$ is invertible if and only
  if its degree-zero term is an invertible element of $A$.

  (ii) It follows from (i) that $1 + yh$ is invertible for all $h \in
  B$.  Hence $y$ is contained in the Jacobson radical $J(B)$.
  Furthermore, it now follows from the natural isomorphism of
  $A$ onto $B/\langle y \rangle$ that $J(B) = J(A) +
  \langle y \rangle$.

  (iii) The $\langle y \rangle$-adic filtration on $B$ is exhaustive,
  separated, and complete. The associated graded ring $\gr B$ is
  isomorphic to $A[z;\alpha]$, and so $B$ is right noetherian, by
  (\ref{filtration}iv), if $A$ is right noetherian.

  (iv) We will refer to the localization of $B$ at the powers of the
  normal element $y$ as the \emph{skew Laurent series ring\/}
  $A[[y^{\pm 1};\alpha]]$.
\end{note}

\begin{note} \label{J-adic} We now briefly summarize some of the
  readily available information on $R$.

(i) Writing
%%%%%%%%%%%%%%%%%%%%%%%%%%%%%%%%%%%%%%%%%%%%%%%%%%%%%%%%%%%%%%%%%%%%%%%%%%%%%%%
  \[ R  = k[[x_1]][[x_2;\tau_2]]\cdots [[x_n;\tau_n]], \]
%%%%%%%%%%%%%%%%%%%%%%%%%%%%%%%%%%%%%%%%%%%%%%%%%%%%%%%%%%%%%%%%%%%%%%%%%%%%%%
  with $\tau_j(x_i) = q_{ji} x_i$ for all $1 \leq i < j \leq n$, we
  see from (\ref{one_variable}ii) that $J = \langle
  x_1,\ldots,x_n\rangle$ is the Jacobson radical of $R$. It follows
  that $J$ is the unique (left or right) primitive ideal of $R$.

  (ii) The $J$-adic filtration on $R$ is exhaustive, separated, and
  complete.  Combining this data with (i), we see that $R$ is a
  complete local ring with unique primitive factor isomorphic to $k$.

(iii) The associated graded ring of $R$ corresponding to the $J$-adic
filtration is isomorphic to the noetherian domain $\qpx$, and it
therefore follows that $R$ is a noetherian domain.

(iv) It follows from \cite[p.~174, 6.~Theorem (3)]{LiVOy} that $R$ is
Auslander regular. It follows from \cite{Wal} that the right Krull dimension,
classical Krull dimension, and global dimension of $R$ are all
equal to $n$.
\end{note}

\begin{note} \label{closed} In view of (\ref{I-adic}ii) and
  (\ref{J-adic}ii, iii), we see that $R$ is zariskian with respect to
  the $J$-adic filtration, since the $J$-adic filtration is complete
  and since the associated graded ring of $R$ corresponding to the
  $J$-adic filtration is noetherian.

  Consequently, the right ideals of $R$ are closed in the $J$-adic
  topology, following (\ref{I-adic}ii).
\end{note}

\begin{note} \label{grlex} Also required in our later analysis
  (particularly in the proofs of (\ref{H-simple}) and
  (\ref{centrally_generated})) is the \emph{graded lexicographic
    ordering\/} on $\Nn$. To review, let $s =
  (s_1,\ldots,s_n)$ and $t = (t_1,\ldots,t_n)$ be $n$-tuples in
  $\Nn$. The \emph{total degree\/} $\vert s \vert$ is the sum $s_1 +
  \cdots + s_n$. Then $s \ltgrlex t$ either when $\vert s \vert <
  \vert t \vert$ or when $\vert s \vert = \vert t \vert$ and $s$
  precedes $t$ in the lexicographic ordering.

  To see the usefulness of this ordering for our purposes, let $N$
  denote an infinite collection of pairwise distinct $n$-tuples in
  $\Nn$. Then we can write $N = \{ j_1,j_2,\ldots \}$ such that
  $j_\ell\ltgrlex j_{\ell+1}$, for all positive integers $\ell$, and
  such that
%%%%%%%%%%%%%%%%%%%%%%%%%%%%%%%%%%%%%%%%%%%%%%%%%%%%%%%%%%%%%%%%%%%%%%%%%%%%%%%
\[ \lim_{\ell \rightarrow \infty} x^{j_\ell} = 0 \]
%%%%%%%%%%%%%%%%%%%%%%%%%%%%%%%%%%%%%%%%%%%%%%%%%%%%%%%%%%%%%%%%%%%%%%%%%%%%%%%
in the $J$-adic topology on $R$.

Also, given a power series $f \in R$, we can choose $s_1,s_2,\ldots
\in \Nn$ and $c_{s_1},c_{s_2},\ldots \in k$ such that
%%%%%%%%%%%%%%%%%%%%%%%%%%%%%%%%%%%%%%%%%%%%%%%%%%%%%%%%%%%%%%%%%%%%%%%%%%%%%%%
\[f = \sum ^\infty_{i=1} c_{s_i}x^{s_i},\]
%%%%%%%%%%%%%%%%%%%%%%%%%%%%%%%%%%%%%%%%%%%%%%%%%%%%%%%%%%%%%%%%%%%%%%%%%%%%%%%
such that $s_\ell \ltgrlex s_{\ell + 1}$ for all positive integers
$\ell$, and such that $c_{s_1} \ne 0$. We will refer to $s_1$ as the
\emph{graded lexicographic degree\/} of $f$.

When the meaning is clear we will simply use ``$<$'' to denote
``$\ltgrlex$.''
\end{note}

The tools developed so far allow the following observation:

\begin{prop} Let $P$ be a prime ideal of $R$. Then $P \cap \qpx$ is prime.
\end{prop}

\begin{proof} Set $Q = P \cap \qpx$. Let $a,b \in \qpx$, and suppose
  that $a(\qpx)b \subseteq Q$. Choose $f \in R$. For each non-negative
  integer $i$, let $f_i$ denote the sum of the monomials in $f$ of
  total degree no greater than $i$. Since $P$ is closed in the
  $J$-adic topology, and since $af_ib \in Q \subseteq P$ for all $i$,
  we see that
%%%%%%%%%%%%%%%%%%%%%%%%%%%%%%%%%%%%%%%%%%%%%%%%%%%%%%%%%%%%%%%%%%%%%%%%%%%%%
\[ afb = \lim_{i\rightarrow \infty} af_ib \in P.\]
%%%%%%%%%%%%%%%%%%%%%%%%%%%%%%%%%%%%%%%%%%%%%%%%%%%%%%%%%%%%%%%%%%%%%%%%%%%%%
Since $f$ was arbitrary, we see that $aRb \subseteq P$. Therefore,
since $P$ is prime, either $a$ or $b$ is contained in $P$. But then
either $a$ or $b$ is contained in $Q = P\cap \qpx$, and so $Q$ is a
prime ideal of $\qpx$.
\end{proof}

\begin{note} \label{laurent} Now let $X$ denote the multiplicatively
  closed subset of $R$ generated by $1$ and the indeterminates
  $x_1,\ldots,x_n$.  Since each $x_i$ is normal, we see that $X$ is an
  Ore set in $R$. For the remainder we will let $L$ denote the
  \emph{$q$-commutative Laurent series ring\/} $\qlsx = \kq[[x^{\pm
    1}_1,\ldots,x^{\pm 1}_n]]$, obtained via the Ore localization of
  $R$ at $X$. (See, e.g., \cite{GooWar} or \cite{McCRob} for details
  on Ore localizations.)

Each Laurent series in $L$ will have the form
%%%%%%%%%%%%%%%%%%%%%%%%%%%%%%%%%%%%%%%%%%%%%%%%%%%%%%%%%%%%%%%%%%%%%%%%%%%%%%
 \[ \sum _{s \in \Zn} c_s x^s ,\]
%%%%%%%%%%%%%%%%%%%%%%%%%%%%%%%%%%%%%%%%%%%%%%%%%%%%%%%%%%%%%%%%%%%%%%%%%%%%%%%
 for $c_s \in k$, for $s = (s_1,\ldots,s_n) \in \Zn$, and for $x^s =
 x^{s_1}\cdots x^{s_n}$, but with $c_s = 0$ for
 $\min\{s_1,\ldots,s_n\} \ll 0$.

 In particular, for $n > 1$ it is \emph{not\/} the case that $L$ can
 be written as an iterated skew Laurent series ring
%%%%%%%%%%%%%%%%%%%%%%%%%%%%%%%%%%%%%%%%%%%%%%%%%%%%%%%%%%%%%%%%%%%%%%%%%%%%%%
\[ k[[x_1^{\pm 1}]][[x_2^{\pm 1};\tau_2]]\cdots[[x_n^{\pm 1};\tau_n]], \]
%%%%%%%%%%%%%%%%%%%%%%%%%%%%%%%%%%%%%%%%%%%%%%%%%%%%%%%%%%%%%%%%%%%%%%%%%%%%%
which is well known to be a division ring.

We always view $R$ as a subalgebra of $L$. Also inside $L$ is the
$q$-commutative Laurent polynomial ring, $\qlx = \kq[x^{\pm 1}_1,\ldots,x^{\pm
  1}_n]$, the localization of $\qpx$ at the products of the
$x_1,\ldots,x_n$.
\end{note}

\section{Prime Ideals in  $q$-Commutative Laurent Series Rings}

We now give a detailed description of the ideals, and in particular
the prime ideals, of $q$-commutative Laurent series rings. Retain the
notation of the preceding section, continuing to let $R = \qpsx$ and
$L = \qlsx$. Henceforth, we will let $Z$ denote the center of $L$.

We begin by considering the ``obvious'' $n$-torus action on $L$.

\begin{note} \label{H-action} Let $H$ denote the algebraic torus
  $(\ktimes)^n$, and set
%%%%%%%%%%%%%%%%%%%%%%%%%%%%%%%%%%%%%%%%%%%%%%%%%%%%%%%%%%%%%%%%%%%%%%%%%%%%%%
\[h(s) = h_1^{s_1}\cdots h_n^{s_n} ,\]
%%%%%%%%%%%%%%%%%%%%%%%%%%%%%%%%%%%%%%%%%%%%%%%%%%%%%%%%%%%%%%%%%%%%%%%%%%%%%%
for $h = (h_1,\ldots,h_n) \in H$ and $s = (s_1,\ldots,s_n) \in
\Zn$. For each $h \in H$, it is not hard to check that the assignment
%%%%%%%%%%%%%%%%%%%%%%%%%%%%%%%%%%%%%%%%%%%%%%%%%%%%%%%%%%%%%%%%%%%%%%%%%%%%%%
\[ h.x^s = h(s)x^s ,\]
%%%%%%%%%%%%%%%%%%%%%%%%%%%%%%%%%%%%%%%%%%%%%%%%%%%%%%%%%%%%%%%%%%%%%%%%%%%%%%
for $s \in \Nn$, extends linearly and continuously to a $k$-algebra
automorphism of $R$. It also is not hard to check, then, that we
thereby obtain an action of $H$ on $R$ by $k$-algebra automorphisms.
In turn, this action by $H$ on $R$ extends to an action by
automorphisms on the localization $L$ of $R$.

Observe that the $H$-actions on $R$ and $L$ induce $H$-actions on
$\spec R$ and $\spec L$. Also note, when $k$ is infinite, that the set
of $H$-eigenvectors in $R$ and the set of monomials in $R$
coincide. The analogous statement for $L$ also holds true.
\end{note}

\begin{note} \label{Gamma-prime} Let $A$ be a ring and let $\G$ be a
  group of automorphisms of $A$. We will refer to ideals of $A$ stable
  under the $\G$-action as \emph{$\G$-ideals}, and we will say that
  $A$ is \emph{$\G$-simple} if the only $\G$-ideals of $A$ are the
  zero ideal and $A$ itself.  We will say that a $\G$-ideal $I$ of $A$
  (other than $A$ itself) is \emph{$\G$-prime\/} if whenever the
  product of two $\G$-ideals of $A$ is contained in $I$, one of these
  $\G$-ideals must itself be contained in $I$. In particular, if $P$
  is a prime ideal of $A$ that is also a $\G$-ideal, then it
  immediately follows that $P$ is also $\G$-prime. If $A$ is right or
  left noetherian, recall that an ideal $I$ of $A$ is a $\G$-ideal
  exactly when $g(I) = I$ for all $g \in \G$.

  Suppose $A$ is right or left noetherian, and suppose that $I$ is a
  $\G$-prime ideal of $A$.  Then $I$ is a semiprime ideal of $A$ and
  is the intersection of a finite $\G$-orbit $I_1,\ldots,I_t$ of prime
  ideals in $A$ all minimal over $I$; see (e.g.) \cite[II.1.10]{BroGoo1}.
\end{note}

\begin{prop} \label{H-simple} Assume that $k$ is infinite. {\rm (i)}
  Every nonzero $H$-ideal of $R$ is generated by monomials. {\rm (ii)}
  $L$ is $H$-simple.
\end{prop}

\begin{proof} (i) Let $I$ be a nonzero $H$-ideal of $R$, and choose
%%%%%%%%%%%%%%%%%%%%%%%%%%%%%%%%%%%%%%%%%%%%%%%%%%%%%%%%%%%%%%%%%%%%%%%%%%%%%%%
\[0  \ne  f  =  \sum_{i \in \Nn} a_i x^i \in I.\]
%%%%%%%%%%%%%%%%%%%%%%%%%%%%%%%%%%%%%%%%%%%%%%%%%%%%%%%%%%%%%%%%%%%%%%%%%%%%%%%
The first step of our proof is to show, as follows, that each $x^i$
appearing nontrivially in $f$ is contained in $I$.

To start, let $\Psi = \{ i \in \Nn \mid a_i \ne 0 \}$. For notational
convenience we will assume that $\Psi$ is infinite; the case when
$\Psi$ is finite can be handled similarly. Equip $\Psi$ with the
graded lexicographic ordering (see (\ref{grlex})), and write $\Psi =
\{ j_1, j_2, j_3, \ldots \}$ such that $j_\ell < j_{\ell+1}$ for all
positive integers $\ell$.
%%%%%%%%%%%%%%%%%%%%%%%%%%%%%%%%%%%%%%%%%%%%%%%%%%%%%%%%%%%%%%%%%%%%%%%%%%%%%%%
So
%%%%%%%%%%%%%%%%%%%%%%%%%%%%%%%%%%%%%%%%%%%%%%%%%%%%%%%%%%%%%%%%%%%%%%%%%%%%%%%
\[ f  =  \sum_{\ell=1}^\infty a_{j_\ell} x^{j_\ell}  .\]
%%%%%%%%%%%%%%%%%%%%%%%%%%%%%%%%%%%%%%%%%%%%%%%%%%%%%%%%%%%%%%%%%%%%%%%%%%%%%%%
Replacing $f$ with $a^{-1}_{j_1}f$, we can assume without loss that
$a_{j_1} = 1$. Set $f_1 = f$ and $a_{j_\ell}(1) = a_{j_\ell}$, for all
$\ell$. Since $k$ is infinite, there exists an $h \in H$ such that
$h(j_1) \ne h(j_2)$. Set $f_2 =$
%%%%%%%%%%%%%%%%%%%%%%%%%%%%%%%%%%%%%%%%%%%%%%%%%%%%%%%%%%%%%%%%%%%%%%%%%%%%%%%
\[\frac{h.f - h(j_2)f}{h(j_1) - h(j_2)}  =  x^{j_1} +
\sum_{\ell = 3}^{\infty}\left(\frac{h(j_\ell) -
    h(j_2)}{h(j_1)-h(j_2)}\right)a_{j_\ell}x^{j_\ell}  =  x^{j_1}
+ \sum_{\ell = 3}^\infty a_{j_\ell}(2)x^{j_\ell}  \in I,\]
%%%%%%%%%%%%%%%%%%%%%%%%%%%%%%%%%%%%%%%%%%%%%%%%%%%%%%%%%%%%%%%%%%%%%%%%%%%%%%%
where
%%%%%%%%%%%%%%%%%%%%%%%%%%%%%%%%%%%%%%%%%%%%%%%%%%%%%%%%%%%%%%%%%%%%%%%%%%%%%%%
\[ a_{j_\ell}(2) = \left(\frac{h(j_\ell) -
    h(j_2)}{h(j_1)-h(j_2)}\right).a_{j_\ell}(1)  .\]
%%%%%%%%%%%%%%%%%%%%%%%%%%%%%%%%%%%%%%%%%%%%%%%%%%%%%%%%%%%%%%%%%%%%%%%%%%%%%%%
Continuing, we get
%%%%%%%%%%%%%%%%%%%%%%%%%%%%%%%%%%%%%%%%%%%%%%%%%%%%%%%%%%%%%%%%%%%%%%%%%%%%%%%
\[ f_m  =  x^{j_1} + \sum_{\ell = m+1}^\infty
a_{j_\ell}(m)x^{j_\ell}  \in I,\]
%%%%%%%%%%%%%%%%%%%%%%%%%%%%%%%%%%%%%%%%%%%%%%%%%%%%%%%%%%%%%%%%%%%%%%%%%%%%%%%
for all integers $m \geq 2$. Note that $f_m \rightarrow
x^{j_1}$ in the $J$-adic topology (again see (\ref{grlex})). By
(\ref{closed}), it then follows that $x^{j_1} \in I$.

Replacing $f$ with $f-x^{j_1} \in I$, we can repeat the above
procedure to show that $x^{j_2} \in I$. Continuing, we see that
$x^{j_1},x^{j_2},x^{j_3},\ldots \in I$. This completes the first step
of our proof.

Now let $M$ be the set of all monic monomials appearing nontrivially
in power series contained in $I$. It follows from the above that $I$
is the smallest closed ideal of $R$ containing $M$. However, by
(\ref{closed}), it then follows that $I$ is the smallest ideal of $R$
containing $M$. In other words, $I$ is generated by $M$, and part (i)
follows.

(ii) Let $I$ be a nonzero $H$-ideal of $L$. Then $I\cap R$ is
  a nonzero $H$-ideal of $R$. Therefore, $I\cap R$ contains a nonzero
  monomial $x^s$ for some $s \in \Nn$. But $x^s$ is invertible in $L$,
  and so $I = L$.
\end{proof}

\begin{cor} \label{cor_to_H-simple} Assume that $k$ is infinite. {\rm
    (i)} Every $H$-orbit of prime ideals in $L$ is (Zariski) dense in
  $\spec L$. {\rm (ii)} $J(L) = 0$.
\end{cor}

\begin{proof} Let $P$ be a prime ideal of $L$, let $O$ be the
  $H$-orbit of $P$ in $\spec L$, and let $I$ be the intersection of
  all of the ideals in $O$. Then $I$ is an $H$-ideal of $L$, and $I$
  is not equal to $L$ itself. Therefore, by (\ref{H-simple}ii), $I =
  0$. It follows that the closure of $O$ in $\spec L$ is $\spec L$,
  and so $O$ is dense in $\spec L$. Part (i) follows. Part (ii) also
  follows, since the intersection of the orbit of a primitive ideal
  must equal the zero ideal.
\end{proof}

\begin{note} Part (ii) of the preceding corollary also follows from
  \cite{Tug}, a study of the Jacobson radical in general skew Laurent
  series rings.
\end{note}

\begin{note} Aside from the last corollary, we are unable at this time
  to provide information on the primitive spectrum of $L$. In
  particular, we ask: Must every primitive ideal of $L$ be maximal?
\end{note}

\begin{note} \label{center} Next, following \cite[\S 1]{GooLet}
  (cf.~\cite[\S 4]{BroGoo2}), we describe $Z$.

  (i) To start, consider the alternating bicharacter
  $\sigma\colon\Zn{\times}\Zn \rightarrow \ktimes$ defined by
%%%%%%%%%%%%%%%%%%%%%%%%%%%%%%%%%%%%%%%%%%%%%%%%%%%%%%%%%%%%%%%%%%%%%%%%%%%%%%
\[ \sigma(s,t)  =  \prod_{i,j = 1}^n \qij^{s_it_j} ,\]
%%%%%%%%%%%%%%%%%%%%%%%%%%%%%%%%%%%%%%%%%%%%%%%%%%%%%%%%%%%%%%%%%%%%%%%%%%%%%
for $s = (s_1,\ldots,s_n), t = (t_1,\ldots,t_n) \in \Zn$.  As noted in
\cite[1.1]{GooLet},
%%%%%%%%%%%%%%%%%%%%%%%%%%%%%%%%%%%%%%%%%%%%%%%%%%%%%%%%%%%%%%%%%%%%%%%%%%%%%
\[ x^s x^t = \sigma(s,t)x^t x^s .\]
%%%%%%%%%%%%%%%%%%%%%%%%%%%%%%%%%%%%%%%%%%%%%%%%%%%%%%%%%%%%%%%%%%%%%%%%%%%%%

(ii) Set
%%%%%%%%%%%%%%%%%%%%%%%%%%%%%%%%%%%%%%%%%%%%%%%%%%%%%%%%%%%%%%%%%%%%%%%%%%%%%
\[ S  =  \{ s \in \Zn  \mid  \text{$\sigma(s,t) = 1$ for
  all $t \in \Zn$} \}.\]
%%%%%%%%%%%%%%%%%%%%%%%%%%%%%%%%%%%%%%%%%%%%%%%%%%%%%%%%%%%%%%%%%%%%%%%%%%%%
As noted in \cite[1.2]{GooLet}, a monomial $x^s$, for $s \in \Zn$, commutes
with all monomials $x^t$, for all $t \in \Zn$, if and only if $s \in S$. It
follows that
%%%%%%%%%%%%%%%%%%%%%%%%%%%%%%%%%%%%%%%%%%%%%%%%%%%%%%%%%%%%%%%%%%%%%%%%%%%%
\[ Z = \left.\left\{ \sum_{s \in S} c_sx^s \in L \right| \text{$c_s
    \in k$, $c_s = 0$ for $\min\{s_1,\ldots,s_n\} \ll 0$} \right\} .\]
%%%%%%%%%%%%%%%%%%%%%%%%%%%%%%%%%%%%%%%%%%%%%%%%%%%%%%%%%%%%%%%%%%%%%%%%%%%%
It will follow from (\ref{second_cor_to_centrally_generated}) that $Z$
is noetherian. However, it is easy to see that $Z$ is a filtered
commutative domain.

(iii) Let
%%%%%%%%%%%%%%%%%%%%%%%%%%%%%%%%%%%%%%%%%%%%%%%%%%%%%%%%%%%%%%%%%%%%%%%%%%%%%
\[ f = \sum_{s \in N} c_s x^s \]
%%%%%%%%%%%%%%%%%%%%%%%%%%%%%%%%%%%%%%%%%%%%%%%%%%%%%%%%%%%%%%%%%%%%%%%%%%%%%
be a Laurent series in $L$, for a suitable index set $N \subseteq
\Zn$. We can write
%%%%%%%%%%%%%%%%%%%%%%%%%%%%%%%%%%%%%%%%%%%%%%%%%%%%%%%%%%%%%%%%%%%%%%%%%%%%%%
\[ f = \sum_{s \in N\cap S} c_s x^s + \sum_{t \in N \setminus S}c_t x^t .\]
%%%%%%%%%%%%%%%%%%%%%%%%%%%%%%%%%%%%%%%%%%%%%%%%%%%%%%%%%%%%%%%%%%%%%%%%%%%%%%
It now follows from (ii) that $Z$ is a right (and left) $Z$-module
direct summand of $L$.

(iv) It is not hard to check that $L$ is finite over $Z$ if and only
if each $\qij$ is a root of unity, if and only if $\rank S = n$.
\end{note}

\begin{note} \label{not_Laurent} Let $b_1,\ldots,b_r$ be a $\Z$-basis
  for $S$.  As noted in \cite[1.3i]{GooLet}, the center of $\qlx$ is
  generated as a $k$-algebra by $(x^{b_1})^{\pm
    1},\ldots,(x^{b_r})^{\pm 1}$ and so is isomorphic to a commutative
  Laurent polynomial ring over $k$ in $r$ many variables. One might
  expect, similarly, that the center of $L$ is a commutative
  Laurent series ring.

  However, the situation for $L$ is more subtle than for $\qlx$. Let
  $\gamma$ be a nonzero, non-root-of-unity in $k$, and consider the
  $q$-commutative Laurent series ring in variables
  $x_1,x_2,x_3$, such that $x_1x_2 = x_2x_1$, $x_1x_3 = \gamma x_3x_1$, and
  $x_2x_3 = \gamma x_3x_2$. We see, in this case, that the central monic
  monomials all have the form $(x_1x_2^{-1})^i$, for integers $i$. It
  is now not hard to see that the center of this $q$-commutative
  Laurent series ring is $k[(x_1x_2^{-1})^{\pm 1}]$ and is
  certainly not a commutative Laurent series ring. Indeed, the
  center in this case is isomorphic to a commutative Laurent
  \emph{polynomial\/} ring. (Our gratitude to Ken Goodearl for this
  observation.)
 \end{note}

\begin{note} \label{change_of_variables} More generally, let
  $\mu_1,\ldots,\mu_n$ be any $\Z$-basis for $\Zn$. Key to the
  approach used throughout \cite{DeCKacPro} and \cite{GooLet} (see
  e.g.~\cite[\S 1.6]{GooLet}) is that the assignment $x^{\mu_i}
  \mapsto y_i$, for $1 \leq i \leq n$, produces an isomorphism from
  $\qlx$ onto $k_r[y_1^{\pm 1},\ldots,y_n^{\pm 1}]$, where $r=
  (r_{ij})$ is the multiplicatively antisymmetric $\nxn$ matrix
  defined by $r_{ij} = \sigma(\mu_i, \mu_j)$.  (Here, $\sigma$ is as
  in (\ref{center}).) In this way, for $\qlx$, one can often reduce to
  the case where $q$ has a particularly desirable form.

  However, changes of variables of this type do not in general extend
  to $L$. To illustrate, consider the commutative Laurent polynomial
  ring $k[X^{\pm 1}]$, in the single variable $X$; note that $X^{-1}
  \mapsto Y$ induces an isomorphism of this algebra onto the Laurent
  polynomial ring $k[Y^{\pm 1}]$, in the single variable $Y$. However,
  this automorphism of Laurent polynomial rings does not extend to an
  automorphism from the Laurent series ring $k[[X^{\pm 1}]]$ onto
  $k[[Y^{\pm 1}]]$.
\end{note}

\begin{note} \label{lemma_history} We now present the ``main lemma''
  of this paper. Analogous conclusions for $q$-commutative Laurent
  polynomial and related rings have appeared in various forms,
  including \cite[5.1]{Bel}, \cite[4.4]{BroGoo2}, \cite[\S
  2]{DeCKacPro}, \cite[2.2]{GooLen}, \cite[1.4]{GooLet},
  \cite[2.1]{Hod}, and \cite[Chapter 11, \S 3]{Pas}.
\end{note}

\begin{lem} \label{centrally_generated} {\rm (i)} Let $I$ be an ideal of
  $L$. Then $I = L(I\cap Z)$.

{\rm (ii)} Let $K$ be an ideal of $Z$. Then $K = (LK)\cap Z$.
\end{lem}

\begin{proof}
(i) If $I = 0$ there is nothing to prove, and so we assume
otherwise. Let $g$ be a nonzero element of $I$, and choose $u \in
\Nn$ such that $x^ug \in I\cap R$. Set $f = x^ug$. Let $S$ be as
in (\ref{center}ii), and fix a transversal $T$ for $S$ in $\Zn$.
By (\ref{center}ii) we can write
%%%%%%%%%%%%%%%%%%%%%%%%%%%%%%%%%%%%%%%%%%%%%%%%%%%%%%%%%%%%%%%%%%%%%%%%%%%%%%
\begin{equation*}
f  =  \sum_{t \in T} x^tz_t ,
\end{equation*}
%%%%%%%%%%%%%%%%%%%%%%%%%%%%%%%%%%%%%%%%%%%%%%%%%%%%%%%%%%%%%%%%%%%%%%%%%%%%%%
where each $z_t$ is a skew Laurent series contained in $Z$, and where
each $x^tz_t \in R$.

Set
%%%%%%%%%%%%%%%%%%%%%%%%%%%%%%%%%%%%%%%%%%%%%%%%%%%%%%%%%%%%%%%%%%%%%%%%%%%%%%
\[T_f = \{t \in T \mid x^tz_t \ne 0 \} .\]
%%%%%%%%%%%%%%%%%%%%%%%%%%%%%%%%%%%%%%%%%%%%%%%%%%%%%%%%%%%%%%%%%%%%%%%%%%%%%%
Then
%%%%%%%%%%%%%%%%%%%%%%%%%%%%%%%%%%%%%%%%%%%%%%%%%%%%%%%%%%%%%%%%%%%%%%%%%%%%%%
\[ f = \sum_{t \in T_f} x^tz_t .\]
%%%%%%%%%%%%%%%%%%%%%%%%%%%%%%%%%%%%%%%%%%%%%%%%%%%%%%%%%%%%%%%%%%%%%%%%%%%%%%
Choose an arbitrary element $t_0$ in the set $T_f$. The main work of
the proof is to prove the following claim: $x^{t_0}z_{t_0} \in I\cap
R$.  Of course, if $f = x^{t_0}z_{t_0}$ there is nothing to prove, and
so we assume otherwise.

To start, we can choose subsets $S_t$ of $S$, for $t \in T_f \setminus
\{t_0\}$, such that
%%%%%%%%%%%%%%%%%%%%%%%%%%%%%%%%%%%%%%%%%%%%%%%%%%%%%%%%%%%%%
\begin{align}
  f \; = \quad x^{t_0}z_{t_0} \, + \sum_{t \in T_f \backslash \{t_0\}}
  x^tz_t \quad &= \quad x^{t_0}z_{t_0} \, + \sum_{t \in T_f \backslash
    \{t_0\}} x^t \sum_{s\in S_t} c_s x^s \notag\\ &= \quad x^{t_0}z_{t_0}
  \, + \sum_{t \in T_f \backslash \{t_0\}} \sum_{s\in S_t} c_s x^t
  x^s, \notag
\end{align}
%%%%%%%%%%%%%%%%%%%%%%%%%%%%%%%%%%%%%%%%%%%%%%%%%%%%%%%%%%%%%
for suitable choices of $c_s\in k$.

Next, regrouping the preceding monomials $c_s x^t x^s$ according to their
total degree, we can write
%%%%%%%%%%%%%%%%%%%%%%%%%%%%%%%%%%%%%%%%%%%%%%%%%%%%%%%%%%%%%%%
\[
f\; = \; x^{t_0}z_{t_0} + \sum_{i\in N} \sum_{t \in T_i} \sum_{s\in S_t^i}
c_s x^t x^s,
\]
%%%%%%%%%%%%%%%%%%%%%%%%%%%%%%%%%%%%%%%%%%%%%%%%%%%%%%%%%%%%%%%%%%%%%%
for sets $N$, $T_i \subseteq T_f \setminus \{t_0\}$, and $S_t^i
\subseteq S$ satisfying the following conditions:
%%%%%%%%%%%%%%%%%%%%%%%%%%%%%%%%%%%%%%%%%
\begin{enumerate}
\item $N = \{1,2,\ldots\}$ is a (possibly finite) index set.

\item $\{ d_i \mid i\in N\}$ is a set of positive integers such that
  $d_i < d_{i+1}$ for all $i \in N$ (if $\vert N \vert > 1$).

\item If $t \in T_i$ and $s \in S_t^i$, then the total degree of the
  monomial $c_sx^tx^s$ is $d_i$.

\end{enumerate}
%%%%%%%%%%%%%%%%%%%%%%%%%%%%%%%%%%%%%%%%%%%%%%%%%%%%%%%%%%%%%%%%%%%%%%%%%%%%
Note that the $T_i$ need not be pairwise disjoint. However,
$T_f\backslash \{t_0\} = \cup T_i$. We now have
%%%%%%%%%%%%%%%%%%%%%%%%%%%%%%%%%%%%%%%%%%%%%%%%%%%%%%%%%%%%%%%%%%%%%%%%%%%%
\[
f=x^{t_0}z_{t_0} + \sum_{i\in N} \sum_{t \in T_i} x^t w_{t, i},
\]
%%%%%%%%%%%%%%%%%%%%%%%%%%%%%%%%%%%%%%%%%%%%%%%%%%%%%%%%%%%%%%%%%%%%%%%%%%%%
where
%%%%%%%%%%%%%%%%%%%%%%%%%%%%%%%%%%%%%%%%%%%%%%%%%%%%%%%%%%%%%%%%%%%%%%%%%%%%
\[w_{t, i}=\sum_{s\in S_t^i} c_s x^s,
\]
%%%%%%%%%%%%%%%%%%%%%%%%%%%%%%%%%%%%%%%%%%%%%%%%%%%%%%%%%%%%%%%%%%%%%%%%%%%%
for $i \in N$ and $t \in T_i$. Note that each $c_s x^s$ is a central
monomial, and so each $w_{t, i}$ is contained in $Z$.

To prove the claim, set
%%%%%%%%%%%%%%%%%%%%%%%%%%%%%%%%%%%%%%%%%%%%%%%%%%%%%%%%%%%%%%%%%%%%%%%%%%%%%
\[ u_1  =  f - x^{t_0}z_{t_0}  =  \sum_{i \in N}\sum_{t \in
  T_i} x^tw_{t, i}\]
%%%%%%%%%%%%%%%%%%%%%%%%%%%%%%%%%%%%%%%%%%%%%%%%%%%%%%%%%%%%%%%%%%%%%%%%%%%%%
and
%%%%%%%%%%%%%%%%%%%%%%%%%%%%%%%%%%%%%%%%%%%%%%%%%%%%%%%%%%%%%%%%%%%%%%%%%%%%%
\[ f_1  = f = x^{t_0}z_{t_0} + u_1  \in I \cap R.\]
%%%%%%%%%%%%%%%%%%%%%%%%%%%%%%%%%%%%%%%%%%%%%%%%%%%%%%%%%%%%%%%%%%%%%%%%%%%%%
Recall that $J = \langle x_1,\ldots,x_n\rangle$ is the augmentation
ideal of $R$ and that each $\sum_{t \in T_i} x^tw_{t, i} $, for $i\in
N$, is a sum of monomials of total degree $d_i$ ($\geq d_1$). Then
$u_1 \in J^{d_1}$. Next, choose some $r \in T_1$. Since $r \ne t_0$, and
since $r$ and $t_0$ are contained in the transversal $T$ for $S$ in
$\Zn$, it follows that $r-t_0 \notin S$. Hence, there exists $v \in
\Zn$ such that
%%%%%%%%%%%%%%%%%%%%%%%%%%%%%%%%%%%%%%%%%%%%%%%%%%%%%%%%%%%%%%%%%%%%%%%%%%%%
\[\sigma(v,t_0) \ne \sigma(v,r),\]
%%%%%%%%%%%%%%%%%%%%%%%%%%%%%%%%%%%%%%%%%%%%%%%%%%%%%%%%%%%%%%%%%%%%%%%%%%%%
where $\sigma$ is as defined in (\ref{center}i). Consider
%%%%%%%%%%%%%%%%%%%%%%%%%%%%%%%%%%%%%%%%%%%%%%%%%%%%%%%%%%%%%%%%%%%%%%%%%%%%
\[\rho_v (f_1) = \frac{x^vf_1x^{-v}-\sigma(v,r)f_1}{\sigma(v,t_0) -
  \sigma(v,r)} .\]
%%%%%%%%%%%%%%%%%%%%%%%%%%%%%%%%%%%%%%%%%%%%%%%%%%%%%%%%%%%%%%%%%%%%%%%%%%%%%
Note that $\rho_v(f_1)$ is contained in $I \cap R$. Also,
%%%%%%%%%%%%%%%%%%%%%%%%%%%%%%%%%%%%%%%%%%%%%%%%%%%%%%%%%%%%%%%%%%%%%%%%%%%%%
\[
\rho_v(f_1) = x^{t_0}z_{t_0} + \sum_{i\in N} \sum_{\substack{t\in
T_i,\\ t\ne r}}\frac{\sigma(v,t) - \sigma(v,r)}{\sigma(v,t_0) -
\sigma(v,r)} x^tw_{t,i}.
\]
%%%%%%%%%%%%%%%%%%%%%%%%%%%%%%%%%%%%%%%%%%%%%%%%%%%%%%%%%%%%%%%%%%%%%%%%%%%%%
Repeating the preceding process, at most $\vert T_1 \vert$-many times,
we obtain
%%%%%%%%%%%%%%%%%%%%%%%%%%%%%%%%%%%%%%%%%%%%%%%%%%%%%%%%%%%%%%%%%%%%%%%%%%%%%
\[ f_2 = x^{t_0}z_{t_0} + u_2 \in I\cap R, \]
%%%%%%%%%%%%%%%%%%%%%%%%%%%%%%%%%%%%%%%%%%%%%%%%%%%%%%%%%%%%%%%%%%%%%%%%%%%%%
with
%%%%%%%%%%%%%%%%%%%%%%%%%%%%%%%%%%%%%%%%%%%%%%%%%%%%%%%%%%%%%%%%%%%%%%%%%%%%%
\[ u_2 = \sum_{\substack{i \in N \\ i \geq 2}}\sum_{t \in T_i
  \setminus T_1} x^tw'_{t,i} ,\]
%%%%%%%%%%%%%%%%%%%%%%%%%%%%%%%%%%%%%%%%%%%%%%%%%%%%%%%%%%%%%%%%%%%%%%%%%%%%%
for suitable $w'_{t,i}\in Z$. Note that each $w'_{t,i}$ appearing in
the preceding will be a nonzero scalar multiple of $w_{t,i}$ and that
$u_2 \in J^{d_2}$. If $u_2 = 0$ then $x^{t_0}z_{t_0} \in I \cap R$,
and the claim follows; so assume otherwise.

Continuing as above, we obtain either $x^{t_0}z_{t_0} \in I \cap R$
(proving the claim) or an infinite sequence
%%%%%%%%%%%%%%%%%%%%%%%%%%%%%%%%%%%%%%%%%%%%%%%%%%%%%%%%%%%%%%%%%%%%%%%%%%%%%
\[ f_i = x^{t_0} z_{t_0} + u_i \in I \cap R ,\quad i=1,2,\ldots \]
%%%%%%%%%%%%%%%%%%%%%%%%%%%%%%%%%%%%%%%%%%%%%%%%%%%%%%%%%%%%%%%%%%%%%%%%%%%%%
converging to $x^{t_0}z_{t_0}$ in the $J$-adic topology. Since the
ideals in $R$ are closed in the $J$-adic topology, as noted in
(\ref{closed}), we see that $x^{t_0}z_{t_0} \in I \cap R$. The claim
follows.

Next, it follows from the claim that $z_{t_0} \in I\cap Z$, and so
%%%%%%%%%%%%%%%%%%%%%%%%%%%%%%%%%%%%%%%%%%%%%%%%%%%%%%%%%%%%%%%%%%%%%%%%%%%%%
\[ x^{t_0}z_{t_0} \in (L(I\cap Z)) \cap R .\]
%%%%%%%%%%%%%%%%%%%%%%%%%%%%%%%%%%%%%%%%%%%%%%%%%%%%%%%%%%%%%%%%%%%%%%%%%%%%%
Recall, however, that $t_0$ was arbitrarily chosen from $T_f$.
Consequently,
%%%%%%%%%%%%%%%%%%%%%%%%%%%%%%%%%%%%%%%%%%%%%%%%%%%%%%%%%%%%%%%%%%%%%%%%%%%%%
\[ x^tz_t \in (L(I\cap Z)) \cap R \]
%%%%%%%%%%%%%%%%%%%%%%%%%%%%%%%%%%%%%%%%%%%%%%%%%%%%%%%%%%%%%%%%%%%%%%%%%%%%%
for all $t \in T_f$.

Following (\ref{grlex}), and recalling that $T$ is a transversal for $S$
in $\Zn$, we can write
%%%%%%%%%%%%%%%%%%%%%%%%%%%%%%%%%%%%%%%%%%%%%%%%%%%%%%%%%%%%%%%%%%%%%%%%%%%%%
\[ \{ x^tz_t \mid t \in T_f \} \]
%%%%%%%%%%%%%%%%%%%%%%%%%%%%%%%%%%%%%%%%%%%%%%%%%%%%%%%%%%%%%%%%%%%%%%%%%%%%%
as a (possibly finite) set
%%%%%%%%%%%%%%%%%%%%%%%%%%%%%%%%%%%%%%%%%%%%%%%%%%%%%%%%%%%%%%%%%%%%%%%%%%%%%
\[ \{ x^{t_1}z_{t_1}, x^{t_2}z_{t_2},\ldots \}, \]
%%%%%%%%%%%%%%%%%%%%%%%%%%%%%%%%%%%%%%%%%%%%%%%%%%%%%%%%%%%%%%%%%%%%%%%%%%%%%
for integers $1,2,\ldots$, with the graded lexicographic degree of
$x^{t_j}z_{t_j}$ less than that of $x^{t_{j'}}z_{t_{j'}}$ whenever $j <
  j'$. Again following (\ref{grlex}), the sequence of partial sums
%%%%%%%%%%%%%%%%%%%%%%%%%%%%%%%%%%%%%%%%%%%%%%%%%%%%%%%%%%%%%%%%%%%%%%%%%%%%%%
\[ x^{t_1}z_{t_1} + \cdots + x^{t_j}z_{t_j} \]
%%%%%%%%%%%%%%%%%%%%%%%%%%%%%%%%%%%%%%%%%%%%%%%%%%%%%%%%%%%%%%%%%%%%%%%%%%%%%%
converges to $f$ in the $J$-adic topology on $R$.  Therefore
%%%%%%%%%%%%%%%%%%%%%%%%%%%%%%%%%%%%%%%%%%%%%%%%%%%%%%%%%%%%%%%%%%%%%%%%%%%%%%
\[f \in L(I\cap Z) \cap R,\]
%%%%%%%%%%%%%%%%%%%%%%%%%%%%%%%%%%%%%%%%%%%%%%%%%%%%%%%%%%%%%%%%%%%%%%%%%%%%%%
since the ideal $L(I\cap Z) \cap R$ is closed in $R$.

Now note that
%%%%%%%%%%%%%%%%%%%%%%%%%%%%%%%%%%%%%%%%%%%%%%%%%%%%%%%%%%%%%%%%%%%%%%%%%%%%%
\[ g  =  x^{-u}f \in L(L(I\cap Z) \cap R) \subseteq L(I \cap Z) .\]
%%%%%%%%%%%%%%%%%%%%%%%%%%%%%%%%%%%%%%%%%%%%%%%%%%%%%%%%%%%%%%%%%%%%%%%%%%%%%
Therefore, $I \subseteq L(I\cap Z)$. Of course $L(I\cap Z) \subseteq
I$, and (i) follows.

(ii) This follows from (\ref{center}iii).
\end{proof}

We now record some applications.

\begin{cor} \label{first_cor_to_centrally_generated} $Z$ is a
  noetherian domain.
\end{cor}

\begin{proof} We know that $L$ is a noetherian domain, and from
  (\ref{centrally_generated}) we know that there is an inclusion
  preserving bijection between the ideals of $Z$ and the ideals of $L$.
\end{proof}

The second application is a partial analogue to
\cite[1.3]{McCPet}.

\begin{cor} \label{second_cor_to_centrally_generated} If $S$ is trivial then
  $L$ is simple.
\end{cor}

\begin{proof} If $S$ is trivial then $Z = k$, and it then follows from
  (\ref{centrally_generated}i) that $L$ is simple.
\end{proof}

The following is analogous to \cite[1.5]{GooLet}.

\begin{prop} \label{homeomorphism} The assignments $P \mapsto P\cap Z$, for $P
  \in \spec L$, and $Q \mapsto LQ$, for $Q \in \spec Z$, provide mutually
  inverse homeomorphisms between $\spec L$ and $\spec Z$.
\end{prop}

\begin{proof} It follows from (\ref{centrally_generated}) that the
  assignments $I \mapsto I\cap Z$, for ideals $I$ of $L$, and $K
  \mapsto LK$, for ideals $K$ of $Z$, provide mutually inverse,
  inclusion preserving bijections between the sets of ideals of $L$
  and $Z$.

  Next, it is a well known elementary fact that a prime ideal in any
  ring intersects with the center at a prime ideal, and so $P\cap Z$
  is a prime ideal of $Z$ for all prime ideals $P$ of $L$. Thus $P
  \mapsto P\cap Z$ produces a map from $\spec L$ to $\spec Z$.

  Now let $Q$ be a prime ideal of $Z$. To show that $LQ$ is a prime
  ideal of $L$, suppose that $I_1$ and $I_2$ are ideals of $L$ such
  that $I_1I_2 \subseteq LQ$.  By (\ref{centrally_generated}i), there
  exist ideals $K_1$ and $K_2$ of $Z$ such that $I_1 = LK_1$ and $I_2
  = LK_2$. So
%%%%%%%%%%%%%%%%%%%%%%%%%%%%%%%%%%%%%%%%%%%%%%%%%%%%%%%%%%%%%%%%%%%%%%%%%%%%%%%
\[ L(K_1K_2) = (LK_1)(LK_2) = I_1I_2 \subseteq LQ .\]
%%%%%%%%%%%%%%%%%%%%%%%%%%%%%%%%%%%%%%%%%%%%%%%%%%%%%%%%%%%%%%%%%%%%%%%%%%%%%%%
Then, by (\ref{centrally_generated}ii),
%%%%%%%%%%%%%%%%%%%%%%%%%%%%%%%%%%%%%%%%%%%%%%%%%%%%%%%%%%%%%%%%%%%%%%%%%%%%%%%
\[ K_1K_2 = (L(K_1K_2))\cap Z \subseteq (LQ)\cap Z = Q.\]
%%%%%%%%%%%%%%%%%%%%%%%%%%%%%%%%%%%%%%%%%%%%%%%%%%%%%%%%%%%%%%%%%%%%%%%%%%%%%%%
Since $Q$ is prime, it follows that either $K_1$ or $K_2$ is contained
in $Q$. Therefore, either $I_1 = LK_1$ or $I_2 = LK_2$ is contained in
$LQ$, and we see that $LQ$ is a prime ideal of $L$. Hence, $Q \mapsto
LQ$ produces a map from $\spec Z$ to $\spec L$.

It now follows from the first paragraph that the assignments $P
\mapsto P\cap Z$, for prime ideals $P$ of $L$, and $Q \mapsto LQ$, for
prime ideals $Q$ of $Z$, provide mutually inverse, Zariski continuous
bijections between $\spec L$ and $\spec Z$. The proposition follows.
\end{proof}

\section{Prime Ideals in $q$-Commutative Power Series Rings}

Our aim now is to systematically develop a detailed description of the
prime spectrum of $R$. Using the results we obtained in the preceding
sections for $q$-commutative Laurent series rings, our approach now
largely follows -- and in many cases mimics -- the studies of
$q$-commutative polynomial rings found in \cite[\S 4]{BroGoo2} and
\cite{GooLet}. Moreover, much of the theory developed in this section
is also analogous to that for various more complicated (finitely
generated) quantum function algebras; see (e.g.)  \cite{BroGoo1} for
details.

\begin{note} \label{stratification} We start by giving an account of
  the ``obvious'' stratification of $\spec R$. Let $W$ be the set of
  subsets of $\{1,\ldots,n\}$.

  (i) For each $w \in W$, let $J_w$ be the ideal of $R$ generated by
  the indeterminates $x_i$ for $i \in w$, and let $R_w = R/J_w$.  Set
  $n_w = n - \vert w\vert$, and observe that $R_w$ is isomorphic
  to a $q$-commutative power series ring (for a replacement of $q$ by
  a suitable $n_w{\times}n_w$ matrix) in $n_w$ many variables. In
  particular, $R_w$ is a domain, and so each $J_w$ is completely
  prime. Next, let
$X_w$ be the multiplicatively closed subset of $R_w$ generated by $1$
and the images of the $x_i$ for $i \notin w$. Set
%%%%%%%%%%%%%%%%%%%%%%%%%%%%%%%%%%%%%%%%%%%%%%%%%%%%%%%%%%%%%%%%%%%%%%%%%%%%%%
\[ L _w  =  R_wX_w^{-1} ,\]
%%%%%%%%%%%%%%%%%%%%%%%%%%%%%%%%%%%%%%%%%%%%%%%%%%%%%%%%%%%%%%%%%%%%%%%%%%%%%%
the Ore localization of $R_w$ at the set $X_w$ (which consists of
normal regular elements of $R_w$). Then $L_w$ is isomorphic to a
$q$-commutative Laurent series ring in $n_w$ many variables (again
using a suitable replacement of $q$).

(ii) For each $w \in W$, let
%%%%%%%%%%%%%%%%%%%%%%%%%%%%%%%%%%%%%%%%%%%%%%%%%%%%%%%%%%%%%%%%%%%%%%%%%%%%%%%
\[ \spec_w R  =  \left\{ P \in \spec R \mid x_i \in P
  \Leftrightarrow i \in w \right\} .\]
%%%%%%%%%%%%%%%%%%%%%%%%%%%%%%%%%%%%%%%%%%%%%%%%%%%%%%%%%%%%%%%%%%%%%%%%%%%%%%%
Then
%%%%%%%%%%%%%%%%%%%%%%%%%%%%%%%%%%%%%%%%%%%%%%%%%%%%%%%%%%%%%%%%%%%%%%%%%%%%%%%
\[ \spec R  =  \bigsqcup _{w \in W} \spec_wR .\]
%%%%%%%%%%%%%%%%%%%%%%%%%%%%%%%%%%%%%%%%%%%%%%%%%%%%%%%%%%%%%%%%%%%%%%%%%%%%%%%
Notice that each $\spec_w R$ is a locally closed subset of $\spec R$;
specifically, $\spec_w R$ is equal to the intersection of the closed
set of prime ideals containing $x_i$, for $i \in w$, with the open set
of prime ideals not containing $x_i$, for $i \notin w$. Also, the
closure in $\spec R$ of each $\spec_wR$ is a union of subsets
$\spec_{w'}R$, for suitable $w' \in W$, and so this partition of $\spec R$
is in fact a stratification.  Equip each $\spec_wR$ with the subspace
topology inherited from $\spec R$.

It follows from \cite[Theorem 10.18]{GooWar} that there is a natural
homeomorphism
%%%%%%%%%%%%%%%%%%%%%%%%%%%%%%%%%%%%%%%%%%%%%%%%%%%%%%%%%%%%%%%%%%%%%%%%%%%%%%%
\[ \spec_w R  \xrightarrow{\; \Theta_w \;}  \spec L _w \]
%%%%%%%%%%%%%%%%%%%%%%%%%%%%%%%%%%%%%%%%%%%%%%%%%%%%%%%%%%%%%%%%%%%%%%%%%%%%%%%
obtained via the assignment
%%%%%%%%%%%%%%%%%%%%%%%%%%%%%%%%%%%%%%%%%%%%%%%%%%%%%%%%%%%%%%%%%%%%%%%%%%%%%%%
\[ P \longmapsto (P/J_w)L_w ,\]
%%%%%%%%%%%%%%%%%%%%%%%%%%%%%%%%%%%%%%%%%%%%%%%%%%%%%%%%%%%%%%%%%%%%%%%%%%%%%%%
for $P \in \spec_wR$. Note that $\Theta_w(J_w) = 0$, for all $w \in W$.

(iii) Let $w \in W$, and let $Z_w$ denote the center of $L_w$. By
(\ref{first_cor_to_centrally_generated}), $Z_w$ is a commutative
noetherian domain, and it follows from (\ref{homeomorphism}) that
$\spec_wR$ is homeomorphic to $\spec Z_w$.
\end{note}

Recall the $H$-actions on $R$ and $\spec R$ discussed in
(\ref{H-action}).

\begin{note} \label{H-equivariance} Let $w \in W$, and assume that $k$
  is infinite.

  (i) Since each $x_i$ is an $H$-eigenvector, it follows that $\spec
  _wR$ is a union of $H$-orbits of prime ideals in $\spec R$, and the
  $H$-action on $\spec R$ restricts to an $H$-action on $\spec_wR$.

  (ii) Again set $n_w = n - \vert w\vert$, and let $H_w$ denote the
  $n_w$-torus $(\ktimes)^{n_w}$. Then $H_w$ acts on $L_w$ as in
  (\ref{H-action}), and so $H_w$ also acts on $\spec L_w$. By
  (\ref{H-simple}ii), each $L_w$ is $H_w$-simple.

  (iii) Letting $H$ act on $\spec_w R$ as in (i), there is an
  ``obvious'' surjection of $H$ onto $H_w$ such that the homeomorphism
  $\Theta_w:\spec_w R \rightarrow \spec L_w$, and its inverse, are
  both $H$-equivariant. In particular, by (ii), each $L_w$ is
  $H$-simple.

  (iv) It follows from (iii) or (\ref{cor_to_H-simple}) that each
  $H$-orbit of prime ideals in $\spec_wR$ is dense in
  $\spec_wR$. (Compare e.g.~with \cite[\S 4.3]{BroGoo2}.)
\end{note}

\begin{prop} \label{Jw} Assume that $k$ is infinite. Then $\{ J_w \mid
  w \in W\}$ is the set of $H$-prime ideals of $R$.
\end{prop}

\begin{proof} To start, we saw in (\ref{stratification}i) that each
  $J_w$ is a (completely) prime ideal of $R$. Moreover, each $J_w$ is
  $H$-stable. It therefore follows from (\ref{Gamma-prime}) that each
  $J_w$ is an $H$-prime ideal of $R$.

  Conversely, let $I$ be an arbitrary $H$-prime ideal of $R$. Then, by
  (\ref{Gamma-prime}), $I = I_1 \cap \cdots \cap I_t$ for some finite
  $H$-orbit $I_1,\ldots,I_t$ of prime ideals of $R$. As noted in
  (\ref{H-equivariance}i), each stratum in the stratification
  described in (\ref{stratification}) is a union of $H$-orbits, and so
  $I_1,\ldots,I_t \in \spec_wR$ for some $w \in W$. Therefore,
%%%%%%%%%%%%%%%%%%%%%%%%%%%%%%%%%%%%%%%%%%%%%%%%%%%%%%%%%%%%%%%%%%%%%%%%%%%%%%%
\[\Theta_w(I_1),\ldots,\Theta_w(I_t)\]
%%%%%%%%%%%%%%%%%%%%%%%%%%%%%%%%%%%%%%%%%%%%%%%%%%%%%%%%%%%%%%%%%%%%%%%%%%%%%%%
is a single $H$-orbit of prime ideals of $L_w$, by
(\ref{H-equivariance}iii). Again following (\ref{H-equivariance}iii),
we see that $L_w$ is $H$-simple, and so
%%%%%%%%%%%%%%%%%%%%%%%%%%%%%%%%%%%%%%%%%%%%%%%%%%%%%%%%%%%%%%%%%%%%%%%%%%%%%%%
\[\Theta_w(I_1) \cap \cdots \cap \Theta_w(I_t) = 0\]
%%%%%%%%%%%%%%%%%%%%%%%%%%%%%%%%%%%%%%%%%%%%%%%%%%%%%%%%%%%%%%%%%%%%%%%%%%%%%%%
in $L_w$. Since $L_w$ is prime, it now follows that
$\Theta_w(I_1),\ldots,\Theta_w(I_t)$ (which comprise a single
$H$-orbit) must all equal $0$. Consequently,
%%%%%%%%%%%%%%%%%%%%%%%%%%%%%%%%%%%%%%%%%%%%%%%%%%%%%%%%%%%%%%%%%%%%%%%%%%%%%%%
\[ I = I_1 = \cdots = I_t = J_w ,\]
%%%%%%%%%%%%%%%%%%%%%%%%%%%%%%%%%%%%%%%%%%%%%%%%%%%%%%%%%%%%%%%%%%%%%%%%%%%%%%%
by (\ref{stratification}ii). The proposition follows.
\end{proof}

\begin{note}
  Next, we consider localization and representation theoretic issues,
  in the sense of \cite{BroWar} and \cite{Jat}. Retaining the notation
  of (\ref{stratification}), recall that the prime spectrum of $R$ is
  \emph{normally separated\/} if for each inclusion of prime ideals
  $P_0 \subsetneq P_1$ in $\spec R$ there exists an element $y \in P_1
  \setminus P_0$ such that $Ry = yR + P_0$. Also, for prime ideals $P$
  and $Q$ of $R$, there is a \emph{(second layer) link\/} $P
  \rightsquigarrow Q$ provided $P\cap Q/PQ$ has a nonzero
  $R$-$R$-bimodule factor that is torsionfree as both a left
  $R/P$-module and a right $R/Q$-module; we can thus view $\spec R$ as
  a directed graph, the connected components of which are termed the
  \emph{cliques\/} in $\spec R$.

  A definition of the \emph{right strong second layer condition\/} can
  be found, for example, in \cite[p.~206]{GooWar}. Normal separation
  implies the strong second layer condition; see, for instance,
  \cite[12.17]{GooWar}. We will use ``$\rank(\; )$'' to denote Goldie
  rank.
\end{note}

Our approach in the following theorem is largely based on \cite{BroGoo2}.

\begin{thm} \label{normal_and_link} {\rm (i)} The prime spectrum of
  $R$ is normally separated, and consequently, $R$ satisfies the
  strong second layer condition. {\rm (ii)} Let $G$ denote the group
  of $k$-algebra automorphisms of $R$ generated by the maps $r \mapsto
  x_i r x^{-1}_i$ for $r \in R$ and $1 \leq i \leq n$. If $P$ and $Q$
  are prime ideals of $R$ such that $P \rightsquigarrow Q$, then
  $\tau(P) = Q$ for some $\tau \in G$. Consequently, if $X$ is a
  clique in $\spec R$, then $\rank(R/P) = \rank(R/P')$ for all $P,P' \in
  X$.
\end{thm}

\begin{proof} (i) Suppose that $P_0 \subsetneq P_1$ are prime ideals
  of $R$, and choose $w$ such that $P_0 \in \spec_wR$. In particular,
  $J_w \subseteq P_0 \subsetneq P_1$. If $P_1 \notin \spec_wR$ then
  there exists some $x_i \in P_1 \setminus P_0$; since $Rx_i = x_iR$
  we see that normal separation holds in this case. Now suppose that
  $P_1 \in \spec_wR$. Set $P'_0 = \Theta_w(P_0)$ and $P'_1 =
  \Theta_w(P_1)$, as in (\ref{stratification}ii). By
  (\ref{stratification}ii), $P'_0 \subsetneq P'_1$. Therefore, by
  (\ref{centrally_generated}i), there exists an element $z$ in the
  center $Z_w$ of $L_w$ such that $z \in P'_1 \setminus
  P'_0$. Moreover, in $L_w$, $z$ is regular modulo $P'_0$. Now recall
  from (\ref{stratification}ii) that $L_w$ is the localization of
  $R_w$ at the set $X_w$ (consisting of normal regular elements of
  $R_w$) and so there exists some $u \in X_w$ such that $uz \in
  R_w$. Note that $uz$ is a normal element of $R_w$, since $u$ is
  normal and $z$ is central. Also, it is easy to see that
%%%%%%%%%%%%%%%%%%%%%%%%%%%%%%%%%%%%%%%%%%%%%%%%%%%%%%%%%%%%%%%%%%%%%%%%%%%%%%%
\[ uz \in (P'_1 \cap R_w) \setminus (P'_0 \cap R_w) .\]
%%%%%%%%%%%%%%%%%%%%%%%%%%%%%%%%%%%%%%%%%%%%%%%%%%%%%%%%%%%%%%%%%%%%%%%%%%%%%%%
But $(P'_0 \cap R_w) = P_0/J_w$ and $(P'_1 \cap R_w) = P_1/J_w$. So
let $y$ be the preimage in $R$ of $uz$. Then $y \in P_1 \setminus
P_0$, and $Ry = yR + P_0$. We can now conclude that the prime spectrum
of $R$ is normally separated. Part (i) follows.

(ii) Suppose that $P$ and $Q$ are prime ideals of $R$ such that $P
\rightsquigarrow Q$. It follows (e.g.) from \cite[12.15]{GooWar}, for
all $1 \leq i \leq n$, that $x_i \in P$ if and only if $x_i \in
Q$. Hence, there exists $w \in W$ such that $P,Q \in \spec_wR$.  In
particular, $J_w \subseteq P \cap Q$.

Set $P' = P/J_w$ and $Q' = Q/J_w$. It follows directly from
\cite[2.7]{BroGoo2} that either $P' \rightsquigarrow Q'$ in $\spec
R_w$ or that $\tau(P) = Q$ for some $\tau \in G$.

So assume that $P' \rightsquigarrow Q'$ in $R_w$. Set $P'' = P'.L_w$
and $Q'' = Q'.L_w$. It is left as an exercise, then, to prove that
$P'' \rightsquigarrow Q''$ in $L_w$. However, since every ideal of
$L_w$ is centrally generated, as proved in
(\ref{centrally_generated}i), it follows that $P'' = Q''$. Hence $P =
Q$. Part (ii) follows.
\end{proof}

A more detailed description of the link structure of
  $\spec R$, following the methods of \cite{BroGoo2}, is left to the
  interested reader.

\begin{note} Next, we combine (\ref{normal_and_link}) with
  \cite{Sta},\cite{War} to obtain a specific application to the
  localization theory of $R$. Recall, when $P$ is a prime ideal of $R$,
  that $\C(P)$ denotes the set of elements of $R$ regular modulo
  $P$. If $X$ is a set of prime ideals of $R$, then $\C(X)$ denotes
  the intersection of all of the $\C(P)$ for $P \in X$, and if $\C(X)$
  is a right Ore set then we use $R_X$ to denote the right Ore
  localization of $R$ at $\C(X)$.

  Now suppose that $X$ is a clique of prime ideals of $R$, and suppose
  that $\C(X)$ is a right Ore set in $R$. Following \cite[\S 7.1]{Jat},
  we say that $X$ is \emph{(right) classically localizable\/} provided
  $R_X$ has the following properties: (1) $R_X/PR_X$ is a simple
  artinian ring for all $P \in X$. (2) Every right primitive ideal of
  $R_X$ has the form $PR_X$ for some $P \in X$. (3) The right
  $R_X$-injective hull of an arbitrary simple right $R_X$-module is
  the union of its socle series.
\end{note}

\begin{cor} \label{classical_cliques} Suppose that $k$ is
  uncountable. Then the cliques in $\spec R$ are classically localizable.
\end{cor}

\begin{proof} By \cite[4.5]{Sta} and \cite[Theorem 8]{War}, cliques
  $X$ for which there is an upper bound on the Goldie ranks modulo
  prime ideals in $X$, in rings with the strong second layer
  condition, are classically localizable. The corollary follows, then,
  from (\ref{normal_and_link}).
\end{proof}

We do not know, in general, whether the cliques in $\spec R$ are
classically localizable.

\begin{catenarity} We obtain a second corollary to
  (\ref{normal_and_link}) by combining it with the results and
  techniques of \cite{ChaWuZha}, \cite{WuZha}, and \cite{YekZha}. (Our
  gratitude to James Zhang both for bringing this corollary to our
  attention and for suggesting to us its proof.)

Recall in a given ring that a chain of prime ideals $P_0 \subsetneq P_1 \subsetneq \cdots \subsetneq P_u$, connecting $P_0$ to $P_u$, is \emph{saturated\/} if for all $1 \leq i \leq u$ there does not exist a prime ideal $P'$ such that $P_{i-1} \subsetneq P' \subsetneq P_i$. Further recall that a ring is \emph{catenary\/} if for each of its comparable pairs of prime ideals $P \subsetneq Q$, every saturated chain of prime ideals connecting $P$ to $Q$ has the same length.
\end{catenarity}

\begin{cor} \label{cat} $R$ is catenary.
\end{cor} 

\begin{proof}  (We require familiarity with the terminology of \cite{ChaWuZha}, \cite{WuZha}, and \cite{YekZha}.) To start, as noted in (\ref{J-adic}iv), $R$ is Auslander regular of global dimension $n$ and so in particular is Auslander Gorenstein.  Therefore, since $R$ is noetherian and local, it follows from \cite[4.3]{ChaWuZha} (cf.~\cite[6.3]{Lev}) that $R$ is AS-Gorenstein (Artin-Schelter Gorenstein). 

  Next, from \cite[Proof of Corollary 0.3, p.~302]{ChaWuZha}, it
  follows that there exists an Auslander (in the sense of \cite[\S 2]{YekZha}), pre-balanced dualizing complex over $R$. Moreover, by \cite[Theorem
  0.1]{ChaWuZha}, this dualizing complex must be $\Cdim$-symmetric. Since $R$ is normally separated, by (\ref{normal_and_link}i), we see that $R$ satisfies all of the hypotheses of \cite[Theorem 6.5]{WuZha}, whose conclusion now ensures that $R$ is catenary.
\end{proof}

\begin{specialcase} \label{special_case} We now consider the following
  situation: Assume that the abelian subgroup $\langle \qij \rangle$
  of $\ktimes$ is free of rank $n(n-1)/2$. (This condition will hold,
  for instance, if the $\qij$, for $1 \leq i < j \leq n$, are
  algebraically independent over $k$.) Of course, in this case $k$ is
  infinite.

  (i) Using (\ref{second_cor_to_centrally_generated}) it is not hard
  to show that each $L_w$, for $w \in W$, is simple. Consequently,
  $\spec R = \{ J_w \mid w \in W\}$. In particular, in this special
  case the prime ideals of $R$ are all completely prime.

  (ii) It follows from (i) and (\ref{normal_and_link}ii) that each
  $J_w$ can only be linked to itself. Hence, by
  (\ref{normal_and_link}i) and \cite[2.5]{BroGoo2}, each $J_w$ is
  linked to itself and only to itself. It then follows from
  (\ref{normal_and_link}i) and (e.g.) \cite[14.20]{GooWar} that each
  $J_w$ is classically localizable.

  (iii) It follows from (i) that $\langle x_1 \rangle, \ldots ,
  \langle x_n \rangle$ are exactly the height-one prime ideals of
  $R$. Moreover, each $\langle x_i \rangle$ is generated by the normal
  element $x_i$ of $R$. In \cite{Cha} (cf.~\cite{ChaJor}), Chatters
  defined a noncommutative noetherian \emph{unique factorization
    domain (UFD) \/} to be a noetherian domain in which every
  height-one prime ideal is completely prime and generated by a normal
  element. Hence $R$ is a UFD in this sense. (Recent results on unique
  factorization in quantum function algebras can be found in
  \cite{LauLen},\cite{LauLenRig}.)

  In \cite{Ven}, the first known examples of noetherian,
  noncommutative, complete, local, unique factorization domains of
  global, Krull, and classical Krull dimension $d > 1$ were presented;
  these examples are completed group algebras of certain uniform
  pro-$p$ groups of rank $2$, and for these algebras $d = 2$. However,
  $R = \qpsx$, in the present (suitably generic) special case, provides examples
  of noncommutative, complete, local, unique factorization domains of
  global, Krull, and classical Krull dimension equal to $n$.
\end{specialcase}

\begin{exam} To see that prime factors of $L$ (and of $R$) need not be
  completely prime (in contrast to the preceding special case),
  consider the situation when $n = 2$, when $k \ne \F_2$, when
  $\varepsilon$ is a primitive $\ell$th root of $1$ in $k$ for some
  integer $\ell \geq 2$, and when
%%%%%%%%%%%%%%%%%%%%%%%%%%%%%%%%%%%%%%%%%%%%%%%%%%%%%%%%%%%%%%%%%%%%%%%%%%%%%%%
  \[ q  =  \bmatrix 1 \hfill & \varepsilon \\ \varepsilon^{-1} &
  1 \endbmatrix.\]
%%%%%%%%%%%%%%%%%%%%%%%%%%%%%%%%%%%%%%%%%%%%%%%%%%%%%%%%%%%%%%%%%%%%%%%%%%%%%%%
  In other words, $L$ now is the $k$-algebra of Laurent series in $x =
  x_1$ and $y=x_2$, subject (only) to the relation $yx = \varepsilon
  xy$. By (\ref{center}ii), the center $Z$ of $L$ is $k[[x^{\pm
    \ell},y^{\pm \ell}]]$. Hence, $L$ is generated as a right (or
  left) module over $Z$ by $\{ x^iy^j \mid 0 \leq i,j \leq
  \ell-1\}$. It follows that the Goldie rank of a prime factor of $L$
  can be no greater than $\ell$.

  Now consider the prime ideal $Q = Z.(x^\ell + y^\ell)$ of $Z$. Then
  $P = LQ$ is a prime ideal of $L$, by (\ref{homeomorphism}). It is
  not hard to check that $(x+y)^i \notin P$, for all $1 \leq i \leq
  \ell -1$. However, it is well known that
%%%%%%%%%%%%%%%%%%%%%%%%%%%%%%%%%%%%%%%%%%%%%%%%%%%%%%%%%%%%%%%%%%%%%%%%%%%%%%%
\[ (x+y)^\ell  =  x^\ell + y^\ell,\]
%%%%%%%%%%%%%%%%%%%%%%%%%%%%%%%%%%%%%%%%%%%%%%%%%%%%%%%%%%%%%%%%%%%%%%%%%%%%%%%
and so $(x+y)^\ell$ is contained in $P$. In other words, $x+y$ is
nilpotent modulo $P$ of index $\ell$. In view of the preceding paragraph, we see that the Goldie rank of $L/P$ must be $\ell$.
\end{exam}

%%%%%%%%%%%%%%%%%%%%%%%%%%%%%%%%%%%-End-Body-%%%%%%%%%%%%%%%%%%%%%%%%%%%%%%%%%%

\end{document}